\documentclass{article}

\newtheorem{theorem}{Theorem}
\newcommand{\reals}{\mbox{\bf R}}

\title{Disintegration of cylindrical measures}

\author{Jan Pachl  \\[-6pt] pachl@acm.org}
\date{December 28, 2002\thanks{Extracted from the author's manuscript dated July 20, 1979.}}

\begin{document}
\maketitle

\begin{abstract}
We show that the existence of disintegration for cylindrical measures follows from
a general disintegration theorem for countably additive measures.
\end{abstract}


\section{Notation}

A {\em probability space\/} $ ( X, \: {\cal A}, \: P ) $ is a nonempty set $X$ together with
a sigma-algebra $\cal A$ on $X$ and a probability $P$ on $\cal A$
(that is, a nonnegative countably additive measure with $PX = 1$).
When $\cal C$ is a set of subsets of a given set,
$\sigma({\cal C})$ is the smallest sigma-algebra that contains $\cal C$.
When $\cal A$ and $\cal B$ are sigma-algebras on $X$ and $Y$, respectively,
$ \sigma ( {\cal A} \otimes {\cal B} ) $ is the smallest sigma-algebra making the
canonical projections $ \pi_X : X \times Y \rightarrow X $
and $ \pi_Y : X \times Y \rightarrow Y $ measurable.

When $ ( X, \: {\cal A}, \: P ) $ and $ ( Y, \: {\cal B}, \: Q ) $ are two probability spaces,
a probability $S$ on $ \sigma ( {\cal A} \otimes {\cal B} ) $ is
a {\em joint probability\/} if $ \pi_X [ S ] = P $ and $ \pi_Y [ S ] = Q $.

A {\em lattice on $X$\/} is a class of subsets of $X$ that
is closed under finite unions and finite intersections.

A class $\cal K$ of sets is {\em semicompact\/} if every countable class ${\cal K}_0 \subseteq \cal K$
such that $\bigcap {\cal K}_0 = \emptyset $ contains a finite class
${\cal K}_{00} \subseteq {\cal K}_0$ such that $\bigcap {\cal K}_{00} = \emptyset $.

Let $ ( X, \: {\cal A}, \: P ) $ be a probability space and $ {\cal K} \subseteq {\cal A} $.
We say that $\cal K$ {\em approximates\/} $P$ if for every $ E \in \cal A$ and $ \varepsilon > 0 $
there is a $K \in \cal K$ such that $ K \subseteq E $ and $ P ( E \setminus K ) < \varepsilon $.

All linear spaces are assumed to be over the field $\reals$ of reals.
When $Y$ and $Z$ are locally convex spaces, ${\cal L} ( Y, Z )$
is the set of continuous linear mappings from $Y$ to $Z$.

When $I$ and $J$ are two index sets such that $J \subseteq I$,
the canonical projection from $\reals^I$ onto $\reals^J$ is denoted $p_{IJ}$.

When $I$ is an infinite index set, let ${\cal C}(I)$ be the set of
all subsets of $\reals^I$ of the form $ p_{IJ}^{-1}(C)$ where $J$ is finite,
$\emptyset \neq J \subseteq I$ and $C$ is a compact subset of
$\reals^J$. Then ${\cal C}(I)_\delta$ is the set of all countable
intersections of sets in ${\cal C} (I)$.

\section{Disintegration theorem}

The following theorem is an immediate consequence of Theorem 3.5 in~\cite{Pachl}.

\begin{theorem}
    \label{th1}
Let $ ( X, \: {\cal A}, \: P ) $ and $ ( Y, \: {\cal B}, \: Q ) $ be two probability spaces,
and let $S$ be a joint probability on $ \sigma ( {\cal A} \otimes {\cal B} ) $.
Let $Q$ be complete,
and let $\cal K$ be a semicompact lattice closed under countable intersections and such that
${\cal A} = \sigma({\cal K})$ and $\cal K$ approximates $P$.

Then there exists a family of probabilities $\{ P_y \}_{y \in Y}$ on $\cal A$ such that

(a) for every $ E \in \cal A$, the function $ y \mapsto P_y E $ is $\cal B$-measurable;

(b) for every $ E \in \cal A $ and $ F \in \cal B $ we have
\[
S ( E \times F ) \: = \: \int_F \; P_y E \:\: {\rm d}Q ( y )  .
\]
\end{theorem}

\section{Application to cylindrical measures}

In his study of non-linear images of cylindrical measures,
Kr\'{e}e~\cite{Kree1} raised the problem of the disintegration of cylindrical measures.
Here we show that the disintegration in the sense of {\em 4.A.b} in~\cite{Kree1} always exists.
However, the disintegration constructed here need not be
continuous or measurable in the sense of {\em 4.A.c-e} in~\cite{Kree1}.
In fact, the remark on page 36 in~\cite{Kree2} seems to imply that there is an example,
due to Schwartz, of a cylindrical measure with no measurable disintegration.

The following theorem improves Proposition (34) in~\cite{Kree1}.
The terminology and notation are as in~\cite{Kree1} and~\cite{Schwartz}.

\begin{theorem}
    \label{th2}
Let $Y$ and $Z$ be two locally convex spaces and let $m$ be a cylindrical probability
on $Y \times Z$ whose canonical image $m_1$ on $Y$ is countably additive.
Extend $m_1$ to the complete probability $Q$ defined on a sigma-algebra $\cal B$.

Then there exists a family $ \{ m_y \}_{y\in Y} $ of cylindrical probabilities on $Z$
such that for every positive integer $n$, every $ b \in {\cal L} ( Z , \reals^n ) $ and
every Borel set $E \subseteq \reals^n $,

(a) the function $ y \mapsto b [ m_y ] E $ is $\cal B$-measurable on $Y$;

(b) for every $F \in \cal B$ we have
\[
( I_Y \times b ) [ m ] ( F \times E ) = \int_{F} b [ m_y ] E \:\: {\rm d} Q ( y ) .
\]
\end{theorem}

The proof of the theorem uses the following result about projective limits of measures.
The original version of this result is due to Bochner.
For a proof of a more general result, see e.g.~\cite{Marczewski}.

\begin{theorem}
    \label{th3}
Let $I$ be an infinite index set. Then \\
(a) ${\cal C} (I)_\delta$ is a semicompact class of subsets of $\reals^I$; \\
(b) if $m$ is a nonnegative finitely additive measure on the algebra generated by ${\cal C} (I)$ and
the image $ p_{IJ} [ m ] $ in $R^J$ is countably additive for each finite $ J \subseteq I $
then $m$ is countably additive and its unique extension to a countably additive measure on the sigma-algebra
$\sigma( { \cal C } (I) )$ is approximated by ${\cal C} (I)_\delta$.
\end{theorem}

{\bf Proof of Theorem~\ref{th2}.}
Schwartz (\cite{Schwartz}, pp. 177-180) points out that cylindrical measures on $Z$ are in one-to-one
correspondence with cylindrical measures on $\reals^I$ for some index set $I$.
Namely, $Z$ is embedded in the algebraic dual $Z'^\ast$ of its topological dual $Z'$, and
$Z'^\ast$ is isomorphic to $\reals^I$ when $I$ is chosen of the same cardinality as an algebraic basis of $Z'$.
This yields an injective continuous linear mapping $h: Z \rightarrow \reals^I$ such that $ \mu \mapsto h[\mu] $
is a one-to-one correspondence between cylindrical probabilities $\mu$ on $Z$ and cylindrical probabilities
on $ \reals^I$.
In the following we choose one such $h$ and keep it fixed.

By Theorem~\ref{th3}, every cylindrical probability on $\reals^I$
is countably additive, and we have a one-to-one correspondence
between cylindrical probabilities on $Z$ and countably additive
probabilities on the sigma-algebra $\sigma( { \cal C } (I) )$ in
$\reals^I$.

Consider the cylindrical probability $ m_0 = ( I_Y \times h ) [ m ] $ on the space $ Y \times \reals^I $
and its projections on $Y$ and $\reals^I$.
The projection $ \pi_Y [ m_0 ] = m_1 $ is countably additive by the assumption.
The projection $ \pi_{\reals^I}$ is countably additive as shown above; let $P$ be its unique
countably additive extension to the sigma-algebra $ {\cal A} = \sigma( { \cal C } (I) )$.

By a result of Marczewski~\cite{Marczewski}, $m_0$ is countably additive.
Thus $m_0$ has a unique extension to a countably additive probability $S$ on the sigma-algebra
$ \sigma ( {\cal A} \otimes {\cal B} ) $.

Hence we can apply Theorem~\ref{th1} with $ X = \reals^I $ and $\cal A$, $P$, $Y$, $\cal B$, $Q$ and $S$ as
defined.
For each probability $P_y$ on $\cal A$ obtained from Theorem~\ref{th1} let $m_y$ be the unique
cylindrical probability on $Z$ such that $ h [ m_y ] = P_y $.

Since $R^I$ identifies with $Z'^\ast$, every $ b \in {\cal L} ( Z , \reals^n ) $
factors through $h$,
and it is now straightforward to verify the properties of $m_y$ stated in the theorem.

{\bf Remark.}
Note that the proof does not use the fact that $Y$ is a vector space.
Thus the same result holds true, with the same proof,
in the more general case of $Y$-cylindrical probabilities (\cite{Kree2}, p.~34),
for a topological or measurable space~$Y$.



\begin{thebibliography}{9}

\bibitem{Kree1}
   {\sc P. Kr\'{e}e.}
   \'{E}quations lin\'{e}aires \'{a} coefficients al\'{e}atoires.
   {\sl Symposia Mathematica Vol. VII\/}, pp. 515-546.
   Academic Press 1971.

\bibitem{Kree2}
   {\sc P. Kr\'{e}e.}
   Accouplement de processus lin\'{e}aires.
   Images de probabilit\'{e}s cylindriques par certaines applications non-lin\'{e}aires.
   {\sl S\'{e}minaire Goulaouic-Schwartz 1971-72.\/}
   Ecole Polytechnique Paris.

\bibitem{Marczewski}
   {\sc E. Marczewski.}
   On compact measures.
   {\sl Fund. Math. 40\/} (1953), 113-124.

\bibitem{Pachl}
   {\sc J. Pachl.}
   Disintegration and compact measures.
   {\sl Math. Scand. 43\/} (1978), 157-168.

\bibitem{Schwartz}
   {\sc L. Schwartz.}
   Radon measures on arbitrary topological spaces and cylindrical measures.
   Oxford University Press 1973.
\end{thebibliography}
\end{document}